\documentclass[12pt]{amsart}
\usepackage{graphicx}
\usepackage{graphics}
\usepackage{amsmath}
\usepackage{amscd}
\usepackage{latexsym}
\usepackage{subfigure}
\usepackage{caption}
\usepackage{hyperref}
\usepackage{placeins}

\begin{document}

\textwidth 5.9in
\textheight 7.9in

\evensidemargin .75in
\oddsidemargin .75in

\newtheorem{Thm}{Theorem}
\newtheorem{Lem}[Thm]{Lemma}
\newtheorem{Cor}[Thm]{Corollary}
\newtheorem{Prop}[Thm]{Proposition}
\newtheorem{Rm}{Remark}

\def\a{{\mathbb a}}
\def\C{{\mathbb C}}
\def\A{{\mathbb A}}
\def\B{{\mathbb B}}
\def\D{{\mathbb D}}
\def\E{{\mathbb E}}
\def\R{{\mathbb R}}
\def\P{{\mathbb P}}
\def\S{{\mathbb S}}
\def\Z{{\mathbb Z}}
\def\O{{\mathbb O}}
\def\H{{\mathbb H}}
\def\V{{\mathbb V}}
\def\Q{{\mathbb Q}}
\def\Cn{${\mathcal C}_n$}
\def\CM{\mathcal M}
\def\CG{\mathcal G}
\def\CH{\mathcal H}
\def\CT{\mathcal T}
\def\CF{\mathcal F}
\def\CA{\mathcal A}
\def\CB{\mathcal B}
\def\CD{\mathcal D}
\def\CP{\mathcal P}
\def\CS{\mathcal S}
\def\CZ{\mathcal Z}
\def\CE{\mathcal E}
\def\CL{\mathcal L}
\def\CV{\mathcal V}
\def\CW{\mathcal W}
\def\IC{\mathbb C}
\def\IF{\mathbb F}
\def\IK{\mathcal K}
\def\IL{\mathcal L}
\def\IP{\bf P}
\def\IR{\mathbb R}
\def\IZ{\mathbb Z}

\title{On an infinite order cork automorphisms}
\author{Selman Akbulut}
\thanks{Partially supported by NSF grants DMS 0905917}
\keywords{}
\address{Department  of Mathematics, Michigan State University,  MI, 48824}
\email{akbulut@msu.edu }
\subjclass{58D27,  58A05, 57R65}
\date{\today}
 
\begin{abstract} 
Here we give a concrete description of the cork automorphism $f:\partial W\to \partial W$ of the infinite order loose-cork $(W,f)$, defined in \cite{a2}. It is obtained by concatenating the defining ribbon disk of $W$ in $B^4$ by an infinite order isotopy of the boundary knot. \end{abstract} 
\maketitle

\setcounter{section}{-1}

\vspace{-.4in}

\section{Introduction }\label{Construction}

\begin{figure}[ht]  \begin{center}
 \includegraphics[width=.3\textwidth]{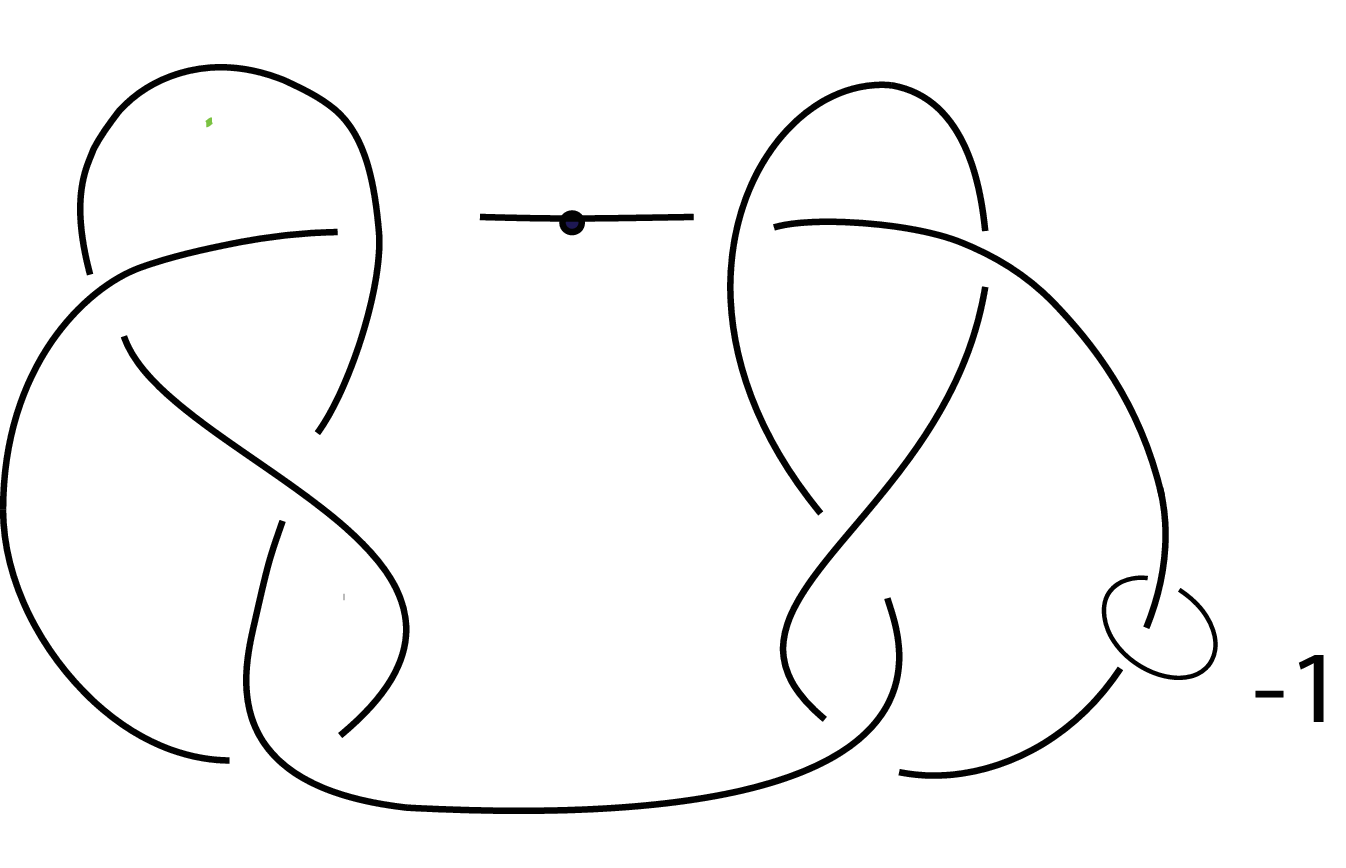}       
\caption{W}      \label{s1} 
\end{center}
 \end{figure}

Recall the infinite order loose-cork $(W,f)$ defined in \cite{a2} and \cite{g}. One definition of $W$ is that, it is the contractible manifold obtained by blowing down $B^{4}$ along the ribbon disk $D\subset B^4$ bounding the  knot $K\#-K$, where $K$ is the figure $8$ knot as in Figure~\ref{s1}  (cf. 6.2 of \cite{a1})

\begin{figure}[ht]  \begin{center}
 \includegraphics[width=.55\textwidth]{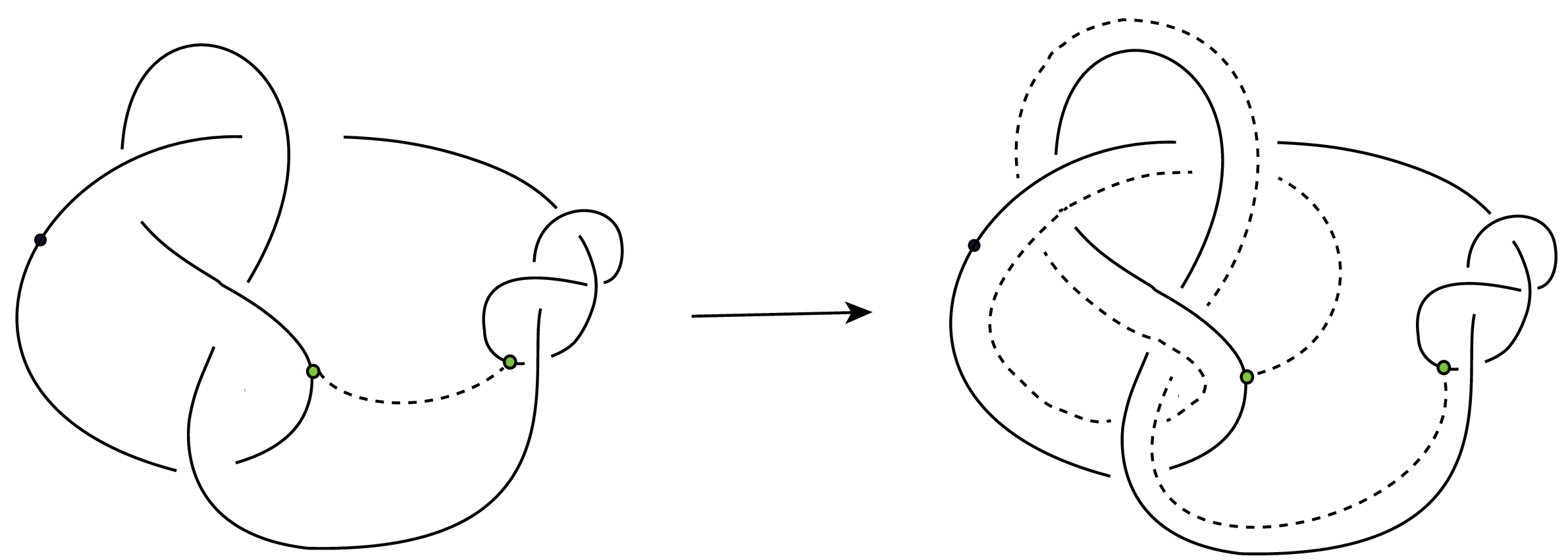}       
\caption{$W_{n}$, for $n=2$}      \label{s3} 
\end{center}
 \end{figure}

What follows from the discussion of \cite{a3} is that, the infinite order cork twisting map  $f_{n}: \partial W \to \partial W$ of this cork is induced from the diffeomorphism  $\partial W \to \partial W_{n}$, where $W_n \approx W$ is the contractible manifold obtained by blowing down $B^{4}$ along the canonical disk $D_{n}\subset B^4$, which is obtained by concatenating $D$ with the concordance of $K\#-K$ along the collar $S^{3}\times [0,1]$ of the boundary, induced by order $n$ ``swallow follow'' isotopy as shown in Figure~\ref{s3}. This isotopy extends to an embient isotopy $F_{t}: B^{3}\times [0,1]\to B^{3}\times [0,1]$ which is  fixed on $B^{3}\times 0\smile S^2\times [0,1]$. 

\begin{Rm}
If the isotopy $F_{t}$ was fixed $B^3 \times 1$, then we would get a nontrivial element of $\pi_{0}$Diff$(B^4, S^3)= \pi_{0}$Diff$(S^4)$, otherwise $(W, f_{n})$ would not be loose-cork.  We would like to thank Rob Kirby noticing that this condition was overlooked in the first posting of this paper. Recall that in  \cite{w}  it was shown $\pi_{k}$Diff$(S^4) \ne \pi_{k}SO(5)$ for some $k>1$.
\end{Rm}


 \begin{figure}[ht]  \begin{center}
 \includegraphics[width=.6\textwidth]{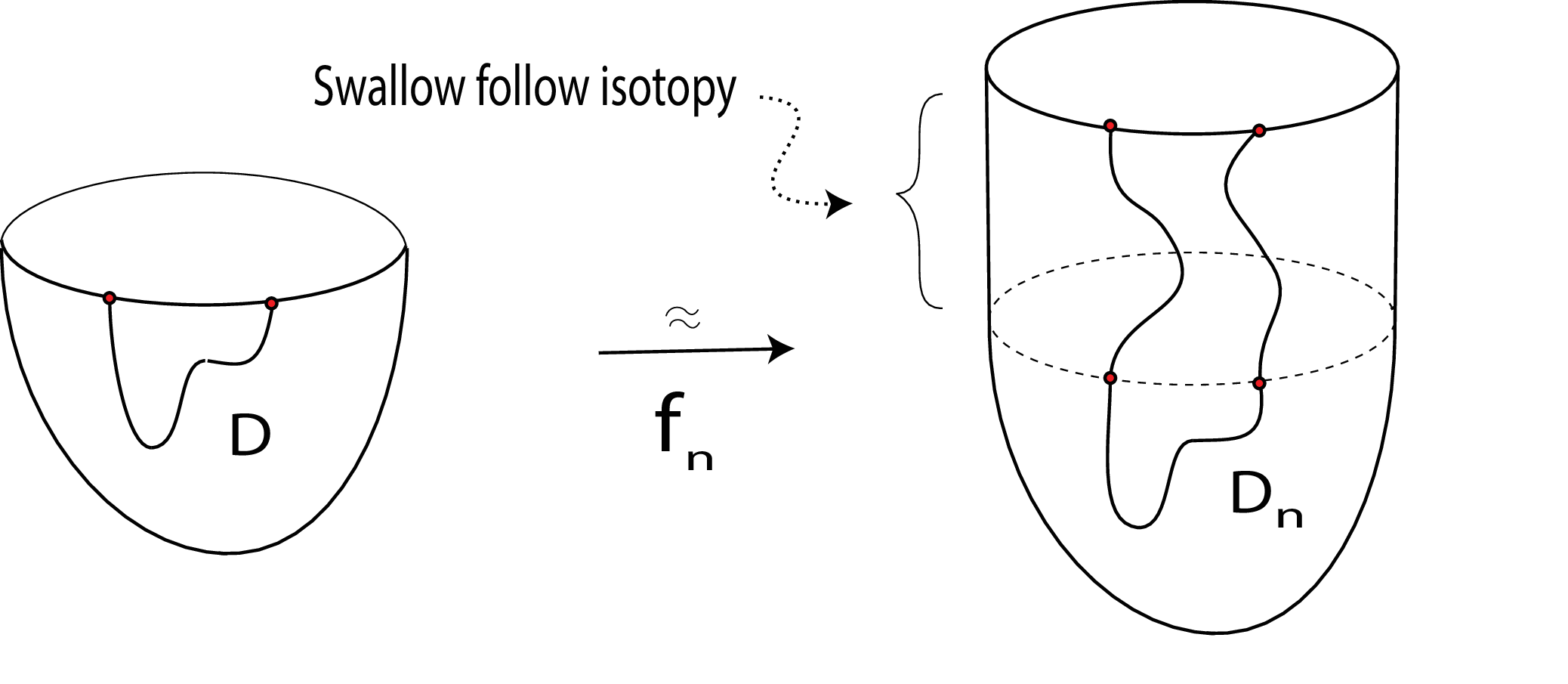}       
\caption{$f: W\to W_{n}$}      \label{s2} 
\end{center}
 \end{figure}

\end{document}